\documentclass[twocolumn,natbib]{svjour3}  
\RequirePackage{fix-cm}
\usepackage[latin2]{inputenc}
\usepackage{hyperref}
\usepackage{times}
\usepackage{t1enc}
\usepackage{amssymb}
\usepackage{amsmath,amsopn}  
\usepackage{t1enc}
\usepackage[mathscr]{eucal}
\usepackage{enumerate}
\usepackage{epsfig,latexsym}
\usepackage{graphicx}
\usepackage{color}
\usepackage{setspace}
\usepackage{marvosym}
\usepackage{epstopdf} 
\usepackage{float}
\journalname{Theoretical and Applied Climatology}

\begin{document}
\title{Applications of threshold models and the weighted bootstrap for Hungarian precipitation data}

\author{  L\'aszl\'o Varga \and P\'al Rakonczai \and Andr\'as Zempl\'eni  }

\institute{L. Varga (\Email ) \and P. Rakonczai \and A. Zempl\'eni \at
              Department of Probability Theory and Statistics, \\   
			  E\"otv\"os Lor\'and University, \\
			  P\'azm\'any P\'eter s\'et\'any 1/C, \\
			  1117 Budapest, Hungary \\
              \email{vargal4@math.elte.hu}           
           \and
           P. Rakonczai \at
              \email{paulo@math.elte.hu}         
           \and
           A. Zempl\'eni \at
			  \email{zempleni@math.elte.hu} }

\date{Received: date / Accepted: date}

\maketitle 

\begin{abstract}  
This paper presents applications of the peaks-over threshold methodology for both the univariate and the recently introduced bivariate case, combined with a novel bootstrap approach. We compare  the proposed bootstrap methods to the more traditional profile likelihood. 

We have investigated 63 years of the European Climate Assessment daily precipitation data for five Hungarian grid points, first separately for the summer and winter months, then aiming at the detection of possible changes by investigating 20 years moving windows. We show that significant changes can be observed both in the univariate and the bivariate cases, the most recent period being the most dangerous, as the return levels here are the highest. We illustrate these effects by bivariate coverage regions.

\end{abstract}
\section{Introduction}

Detection of signs for climate changes are in the focus of recent climatology. There is an abundance of publications in the area of temperature changes. Precipitation is equally important, if we consider its economic effect -- here the extremes play an especially important role, since these are in close connection with dangerous floods or extreme draught. That is why we focus on the extreme value models, for the upper tails of the distribution, i.e. on the estimation of the return levels for extreme daily precipitation. There are some papers which deal with extreme precipitation: \cite{bp}shows that there is an increase in days with extreme precipitation in the region (here the investigted period lasts till 2001). A more recent paper is \cite{ bp2}, where also simulations for the 21st century are analysed. These analyses are based on indices, like annual number of days with precipitation over 10 mm, and do not give estimates for the high quantiles of the underlying distribution.  \cite{hb} found a significant increase in winter precipitation while a decrease in summer values for the Rhine basin in Germany. Besides the indices they also investigate the 90\% quantiles.  Linear trend functions were also fitted in several papers to the annual maximum precipitation, for example see  \cite{kk}. These trends were in most cases not significant (in contrast to trends in temperature indices) and did not show any spatial pattern. 

Our purpose is to show that indeed the spatial patterns do change. This is a question which has not yet been investigated  in the context of precipitation. Mathematically this means the investigation of the joint distribution of the extremes. However, in order to  be able to tackle this question, we need to model the univariate distributions first. The used mathematical model here is the peaks over threshold model for the univariate extremes, as the most common and also theoretically sound model. See for example \cite{dsa}, where several potential models were compared and this peaks-over-threshold model turned out to be the best. 
We also investigate the dependencies among the extreme observations, for this purpose the bivariate peaks over threshold 
models are used. We give details of these models in Section \ref{models}. 

We intend to quantify the uncertainties of our estimators. As we'd need much more observations for the classical limit theorems to give accurate approximations in the extreme value models, we had to find other methods for this purpose.
The bootstrap is a controversial tool in the area of extremes. There are reports on its too short confidence intervals and low coverage probability, see for example \citet{kys}. However, one needs methods for estimating the uncertanties of the estimators and the asymptotic normality is not always a reasonable approximation either. 
We propose an alternative use of the bootstrap resampling technique, which helps to overcome the known problems. It is combined with the profile likelihood method, with quite satistactory results.

\begin{figure} 
\vspace{0mm}
\includegraphics[width=80mm]{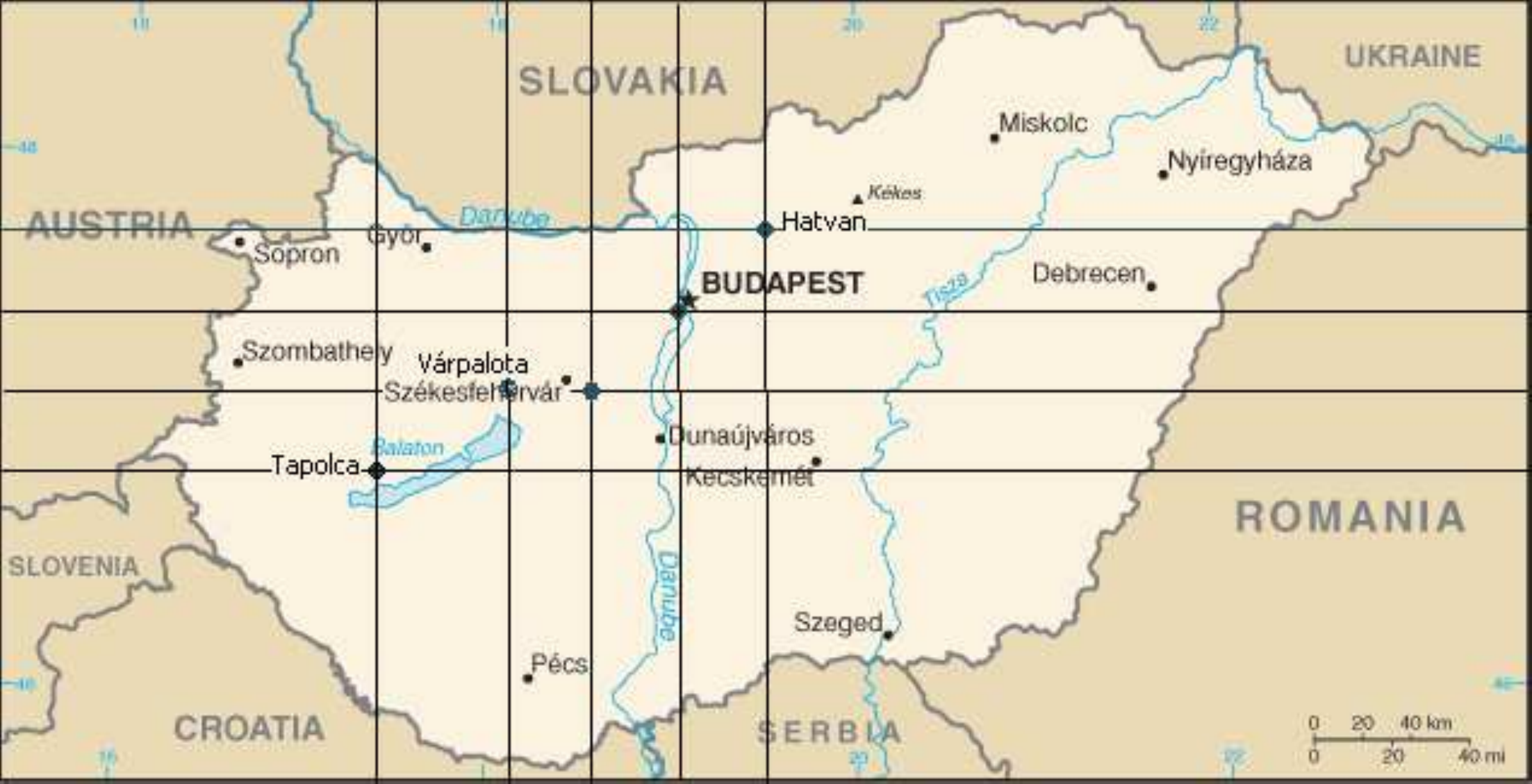}
\vspace{0mm}
\caption{Map of  Hungary with the used grid points (blue circles) and some major cities}
\label{fig:map}
\end{figure}

The observations are 63 years of daily precipitation data for 5 Hungarian grid points, corresponding to the following cities: Tapolca, V\'arpalota, Sz\'ekesfeh\'erv\'ar, Budapest and  Hatvan, from the grid E-OBS ({\url {http://www.ecad.eu}}).  This is a gridded data base, which has been used extensively for climate analysis, see \cite{hhgk}. The quality has been evaluated in \cite{hhnj}, and the results show that it may be considered reliable for most of Central Europe.  Figure \ref{fig:map} depicts the used five grid points together with some major cities on the Hungarian map.  Based on our work we hope to reveal  important characteristics of the joint behaviour of the precipitation data as well as a more exact evaluation of the significance of the trends in the univariate data sets.

In Section \ref{models} first the univariate, then the multivariate threshold models are introduced. Section \ref{boot} contains the bootstrap approach and its applications for extreme value models. In Section \ref{app}  we show the univariate and bivariate applications of the models. Section \ref{conclu} contains the conclusions.


\section{Peaks-over-threshold models}
\label{models}

There are two widely investigated types of approaches to model extreme values. The "classical" models are based on the block maxima of the data. The other, more recent approach focuses on the observations exceeding a certain, high threshold. The latter class of models is called peaks-over-threshold  (POT) models. In this paper we will concentrate on such models, because they allow for the use of more data and turned out to be powerful in quite a few practical problems.  \cite{bv} have used this model to build a spatial pattern for extreme precipitation hazard, based on daily data. Our approach differs from this, as our intention is to determine temporal trends and their significance.

\subsection{The univariate case}
\label{Ch21}
Univariate threshold models have a long history since their introduction 
in the 1970s \citep[see][]{bh, pi}.
Under fairly general regularity conditions the threshold exceedances have an asymptotic distribution. To be more specific,the conditional excess distribution function converges:  
\begin{equation*}
\label{for:exc}
F_u(z)=P(X-u\leq z| X>u)  \overset{u \to \infty}{\longrightarrow}  H(z),
\end{equation*}
where $H(z)$ is a distribution function with real valued parameters $\xi$ (shape) and $\sigma>0$ (scale).
\begin{equation}
\label{for:gpd}
H(z)=
\begin{cases}
1-\left( 1+\frac{\xi z}{\sigma} \right)^{-\frac{1}{\xi}} & \text{if } \xi\neq 0 \\
1-e^{-\frac{z}{\sigma}} & \text{if } \xi=0
\end{cases}.
\end{equation} 
The family defined in (\ref{for:gpd}) is called generalized Pareto distribution (GPD). 
Depending on the parameter $\xi$, this distribution includes three types of distribution familes, 
which is summarized in Table \ref{tab:gpd}. \\
\begin{table}[ht]
\caption{Distribution families covered by the GPD distribution}
\label{tab:gpd}     
\begin{center} \begin{tabular}{lll}
\hline
& Distribution family & Support \\ \hline
 $\xi>0$ & Ordinary Pareto & $0 \leq z$ \\
 $\xi<0$ & Type II Pareto & $0 \leq z \leq -\frac{\sigma}{\xi}$ \\
 $\xi=0$ & Exponential & $0 \leq z$ \\  \hline
\end{tabular} \end{center}
\end{table}
In statistics, if  $\mathbf{X}_n=(X_1,...,X_n)$ is a sequence of independent observations with identical (unknown) distribution function $F$ (for which we shall use the common abbreviation of i.i.d. sample), then the exceedances will be considered as a sample from the GPD $H$.
In our example, while the original daily observations are clearly dependent, this dependence is practically negligible if only the exceedances are considered. And it is worth mentioning that the theory implying the limit is valid for weak dependence as well.

These distributions have proved to be a suitable model for precipitation data, see for example \cite{dsa}, where several reasonable families were compared and the Pareto distribution was clearly the best. We shall see that in our case both ordinary and type II Pareto distributions appear, showing the difference between places and the seasons.

If $\xi\neq 0$, then the density function of the GPD distribution is
\begin{equation*}
\label{for:dens}
h_{\xi,\sigma}(z)= \frac{1}{\sigma} \left( 1+\frac{\xi z}{\sigma} \right)^{-\frac{1}{\xi}-1},
\end{equation*}
and so we get the log-likelihood function 
\begin{equation}
\label{for:llh}
l(\xi,\sigma|\mathbf{X}_n)= \sum\limits_{i=1}^n \log h_{\xi,\sigma}(X_i).
\end{equation}
We denote the maximum likelihood (ML) estimate of the parameters by $\hat{\sigma},\hat{\xi}$.

In meteorology return levels corresponding to certain return periods -- for instance 10, 20 or 50 years -- are especially important. That's why the log-likelihood function is often parameterized in terms of $\xi$ and the quantile 
function $H^{-1}$ \footnote{The (old) parameter $\sigma$ in terms of $\xi$ and $H^{-1}(q)$: \linebreak
 $\sigma= \sigma(\xi,H^{-1}(q)) =\frac{\xi H^{-1}(q)}{(1-q)^{-1}-1}.$}:
\begin{flalign*}
\label{for:llh2}
 l(\xi,H^{-1}(q)|\mathbf{X}_n) 
 = \sum\limits_{i=1}^n \log h_{\xi,H^{-1}(q)}(X_i), \text{\quad where}
\end{flalign*}
{\small
\begin{equation*}
h_{\xi,H^{-1}(q)}(z) = 
\frac{(1-q)^{-1}-1}{\xi H^{-1}(q)} \left( 1+z \frac{(1-q)^{-1}-1}{H^{-1}(q)} \right)^{-\frac{1}{\xi}-1}.
\end{equation*}
}
Now the ML estimates are denoted by $\hat{\xi}$ and $\widehat{H^{-1}(q)}$. In the applications we are going to apply this parameterization. Let us note that the ML estimators have the usual rate of convergence and normal limit if $\xi>-0.5$, which is the case in almost all applications (including ours).

 It is important to note that in general 
if we have $n$ observations over the threshold in $l$ 
years\footnote{This means in average $\frac{n}{l}$ observations in a single year.}
and the return period of interest is $m$ years then the corresponding quantile is
 $q=1-\frac{1}{m} \frac{l}{n}.$
For example if there are 10000 such observations in 100 years then the 50 years return level is the 0.9998-quantile of the distribution, which models the single exceedances.

\subsection{The bivariate case}
\label{Ch22}
The multivariate counterpart of the peaks-over-threshold models has been recently developed.  For case of simplicity we introduce the bivariate case only -- we shall apply this model for the data. Bivariate threshold models can be defined in two different ways.

If we claim exceedances in both coordinates (called {\color{white} i} \mbox{BGPD I} model), we usually get simpler models, with nice properties (margi\-nals are univariate GPD etc.), but we may use less data.

If we use all data that exceed the threshold ${\mathbf u}$ in at least one coordinate, we get the so called BGPD II model of \citet{rt}.

This approach can be formulated as follows.
Let $\mathbf{Y}=(Y_1,Y_2)$ denote a random vector, $\mathbf{u}=(u_1,u_2)$ be a suitably high threshold vector and $\mathbf{X} = \mathbf{Y}-\mathbf{u} = (Y_1-u_1,Y_2-u_2)$ be the vector of exceedances. Then the bivariate generalized Pareto distribution (BGPD) for the exceedances $\mathbf{X}$  can be defined by a bivariate extreme value distribution $G$ with non-degenerate margins as
$$
H({\mathbf x}) =\frac{1}{\log G\bigl(0,0 \bigr)} \log \frac{G\bigl(x_1,x_2 \bigr)}{G\bigl(x_1 \wedge 0,x_2 \wedge 0 \bigr)} ,
\label{mgpd} $$
where $0<G(0,0)<1$. 
 
{\bf Remarks}
\begin{enumerate}
\item
 Note that 
the margins $X_i$ are one dimensional GPDs, only under the condition $X_i>0$.
\item All margins are dependent on all
parameters, as the constant factor  $1/ \log G(0, 0)$ remains in the formula. 
So the parameters cannot be interpreted individually.
  \item Some models put weight to the boundaries, so they will not remain absolutely continuous. 
\item The most important advantage of this approach is that we can use more data, which hopefully helps in model fitting.
\end{enumerate}
  
The most popular parametric model
is the (symmetric) logistic model:  \\
$l(v_1,v_2)=(v_1^{1/\alpha}+v_2^{1/\alpha})^{\alpha}, \qquad v_j\geq 0,$ with parameter
$0<\alpha\leq 1$ where $l=-\log G$. Independence occurs when $\alpha =1$, complete dependence for $\alpha \downarrow 0$. It is symmetric, absolutely continuous.

There are of course many more classical as well as new models (see e.g. \cite{rz} for an over\-view).   
\section{Bootstrap}
\label{boot}
The bootstrap is a resampling method developed in the last two decades of the previous century. The main goal of the bootstrap is to extract the maximum amount of information from the data on hand. The basic, i.i.d. bootstrap idea is to produce new samples from the original one via resampling with replacement. The bootstrap is appropriate for estimating the uncertainty of the statistics of interest, computing confidence intervals and $p$-values. There are a lot of different bootstrap techniques, for further details see for example the books \citet{ch} or \citet{dh}. There are also interesting meteorological applications to be found in the literature, for example in a recent paper \cite{us}
 use a spatial variant of the bootstrap methodology for estimating the distribution of annual maximum precipitation.

In the following subsection we will focus on the so-called weighted bootstrap.

\subsection{Weighted bootstrap}
\label{boot1}
The weighted bootstrap -- sometimes called multiplier bootstrap -- is an extension of the bootstrap scheme above.
In this framework the multiplicity of the elements of the original sample in the bootstrap sample are considered to be random variables. These random variables are called bootstrap weights and they will be denoted by $\tau_{ni}$ $(i=1,2,...,n)$, where $n$ is the sample size. So the bootstrap sample -- for integer-valued $\tau$ -- becomes
\begin{equation}
\mathbf{X}^*_n=(\underbrace{X_1,\dots X_1}_{\tau_{n1} \textrm{times}},\dots,\underbrace{X_n,\dots X_n}_{\tau_{nn} \textrm{times}}) .
\end{equation}
The weighted bootstrap can be found in many recent applications:  \citet{bd} approximated the empirical copula process with it; \citet{df}, \citet{fhh} and \citet{rmi} exploited weighted bootstrap in regression analysis.

In case of general $\tau$, the weights are simply applied to the log-likelihood function. In this approach the elements of the log-likelihood function are multiplied by the bootstrap weights the following way:
\begin{flalign*}
\label{for:llh2}
 l^*(\xi,H^{-1}(q)|\mathbf{X}_n) 
 = \sum\limits_{i=1}^n \tau_{ni} \log h_{\xi,H^{-1}(q)}(X_i).
\end{flalign*}
The main theoretical concepts, properties and consequences of this "weighted bootstrap-likelihood" were investigated in \citet{wa}. Recently it has been used in autoregressive models (\citet{bb}) and ARCH processes (\citet{vz}). This likelihood-based approach reduces the computing burden of this sophisticated bootstrap method, which becomes especially important when fitting rather complicated models, like the BGPD II, see subsection \ref{Ch22}. 

Let us suppose the following assumptions \textbf{B1-B5} for the bootstrap weights: 
\begin{itemize}
\item[\textbf{B1}] The weights are independent from the data-generating process.
\item[\textbf{B2}] The weights are nonnegative: $P(\tau _{ni} \geq 0)=1;  \quad  i=1,...,n; \ \ n=1,2,... $
\item[\textbf{B3}] The first two moments of $\tau_{n1}$, $\dots$, $\tau_{nn}$ are finite and equal for any fixed $n$.
\item[\textbf{B4}] $\underset{n \to \infty}{\lim} E \tau _{ni} = 1 \quad i=1,2,... $
\item[\textbf{B5}] $\gamma:=\underset{n \to \infty}{\lim} E \tau_{ni}^2 < \infty \quad i=1,2,... $
\end{itemize}

Several distributions fulfill the five assumptions above, we shall use the multinomial and i.i.d. exponential weights in this paper (it is worth mentioning that the multinomial case corresponds to the classical bootstrap setup). To be more specific, we either use
\begin{equation*}
\label{wei:mult}
(\tau_{n1},...,\tau_{nn}) \sim \text{ Multinom} \left( n;\frac{1}{n},...,\frac{1}{n} \right),
\end{equation*} or
\begin{equation*}
\label{wei:exp}
(\tau_{n1},...,\tau_{nn}) \sim \text{ i.i.d. Exp(1)}.
\end{equation*}

\subsection {Bootstrapping the extremes and profile likelihood}
\label{boot2}
There are different approaches for bootstrapping in the extreme value models: parametric and nonparametric bootstrap, parametric
being the most commonly used method (see \citet{kys} for a survey in the meteorological context). However, our new proposed method -- which is nonparametric -- seems to be much more flexible than
the original one, and as -- for example in our bivariate investigations --the model is not unique, the extremes are not especially
heavy tailed and there is a relatively long data set available, this nonparametric approach is more appealing. 

In meteorology the profile likelihood method is widely applied to construct confidence intervals for return levels (quantiles) or other relevant parameters, see e.g. \cite{du}.
 The used main concept is the so-called \textit{profile} log-likelihood function \citep[see][p. 33-36]{co}
which is 
\begin{equation}
\label{for:pllh}
l_p(H^{-1}(q)|\mathbf{X}_n)= \underset{\xi}{\max} \ \ l(\xi,H^{-1}(q)|\mathbf{X}_n)
\end{equation}
The profile log-likelihood function is the maximized log-likelihood with respect to $\xi$. So it gives the local maxima of the log-likelihood function for different $H^{-1}(q)$ values.

If $\hat{\xi}$ and $\widehat{H^{-1}(q)}$ denotes the ML-estimate of the parameters then under some regularity conditions it is well 
known\footnote{Special case of Theorem 2.6 in \citet{co}.} that
\begin{equation}
\label{for:limdist}
2\left[ l(\hat{\xi},\widehat{H^{-1}(q)}|\mathbf{X}_n)-l_p(H^{-1}(q)|\mathbf{X}_n) \right] 
\end{equation}
has a $\chi^2$ limit distribution with 1 as the degree of freedom as $n\to\infty$. This result can be used to construct confidence intervals for the return levels. The profile intervals are rarely symmetric and in practice they seem to perform really well compared to other methods' intervals.
   
We combined the weighted bootstrap with the profile likelihood method to construct confidence regions for the return levels of the precipitation data on hand. 

The bootstrap version of the profile log-likelihood function is simply the appropriate modification of (\ref{for:pllh}).

Let $\gamma$ be the limit of the second moments of the weights defined in (\textbf{B5}) It can be easily proved that under conditions B1-B5 
\begin{equation}
\label{for:limdistboot}
\frac{2}{\gamma}\left[ l^*(\hat{\xi},\widehat{H^{-1}(q)}|\mathbf{X}_n)-l_p^*(H^{-1}(q)|\mathbf{X}_n) \right],
\end{equation}
also has the same limit distribution as  (\ref{for:limdist}). This asymptotic result can be used to construct confidence regions for the return levels. In the sequel the level of confidence will be denoted by $1-\alpha$, typically chosen as 0.95 or 0.99; the $(1-\alpha)$-quantile of the $\chi^2_1$-distribution by $c_{1-\alpha}$; and the empirical sample by $\mathbf{x}=(x_1,...,x_n)$. So using 
(\ref{for:limdistboot}), the constructed profile bootstrap confidence interval $I^*_{\alpha}$ becomes
\begin{equation}
\begin{split}
I^*_{\alpha}  = &\left\{  H^{-1}(q):  {\color{white} \frac{c}{2}}\right. \\
& \left.  l^*_p(H^{-1}(p)|\mathbf{x}) \geq l^*(\hat{\xi},\widehat{H^{-1}(p)}|\mathbf{x}) -
\frac{\gamma \cdot c_{1-\alpha}}{2} \right\}.
\end{split}
\end{equation}

This confidence region is usually wider than the conventional profile likelihood confidence interval and as we will see, it outperforms the traditional methods applied to our precipitation dataset, especially if we don't confine our scope to the seasons. See subsections \ref{appuni} and \ref{appbi} for the applications.

\section{Applications}
\label{app}

Having shown the tools, now we are in position of actually carrying out the data analysis.

We have used 63 years (1950 to 2012) of daily precipitation data of 5 Hungarian grid points, corresponding to the following cities: Tapolca, V\'arpalota, Sz\'ekesfeh\'erv\'ar, Bu\-da\-pest and  Hatvan, as shown on Figure \ref{fig:map}. The choice was motivated by the
fact that in the bivariate analysis we were interested in the analysis of the dependencies between the daily observations, which turned out to be rather small if we consider sites which are far away. Budapest was a natural choice and the other grid points were chosen nearby, with emphasis to the slightly wetter Transdanubian part of Hungary. We intended to give detailed analysis, that is why we stuck to this limited number of grid points, which was just enough to show that our methods work for different places as well. In quite a few steps of the analysis we do not show all the results, but illustrate the situation by one or two typical figure. Unless otherwise stated this means that the other stations showed similar patterns. The main question we had in mind was if there was an observable and statististically significant change in the precipitation over this period. We did not follow the linear regression approach, proposed by some authors (see \citet{cc} for a recent work in the area), as there may be other forms of changes, which may be observed by nonparametric methods. Another important question we had in mind was the quantification of the uncertainty of our estimators. This is  an important question, as the standard assumptions like asymptotic normality does not hold in the case of the extremes, we are interested in.
 
First we investigated the time series, which shows some seasonality. This is markedly true for the extremes as well.
See Table \ref{tab:arany} for the monthly frequencies of observations beyond the chosen threshold of 10 mm. 

\begin{table*}[ht]
\centering
\caption{Monthly frequencies of observations beyond the chosen threshold of $u$=10 mm (cumulated exceedances)}
\label{tab:arany}     
\begin{tabular}{lccccc}
\hline
Months &       Tapolca & V\'arpalota & Sz\'ekesfeh\'erv\'ar &  Budapest &  Hatvan \\ \hline
January &      33  &           30   &   35   &    37  &      32 \\ 
February &     30  &           31   &   28   &    36  &     31 \\ 
March &        37  &           30   &   15   &    15  &     19 \\ 
April &        60  &           53   &   41   &    46  &    35 \\ 
May &          96  &            80   &   74   &   78  &      79 \\ 
June &         134 &          103  &   104  &    110 &     89 \\ 
July &         124 &            97   &   81   &    79  &     79 \\ 
August &       135 &            105  &   91   &    70  &       75 \\
September &    113 &         91   &   69   &    65  &       71 \\ 
October &      89  &            80   &   67   &    62  &     58 \\ 
November &     102 &         102  &   82   &   91 &   82 \\ 
December &     66  &         56   &   56   &   59  &     49 \\ \hline    
\end{tabular}
\end{table*}

So we decided to investigate first the seasonal models, where especially the summer months (JJA) are the most interesting as convective events may cause heavy rains, but with different extremal characteristics and spatial dependence structure as the typical frontal systems of the other seasons.

\begin{figure}[!ht]
\includegraphics[scale=0.35]{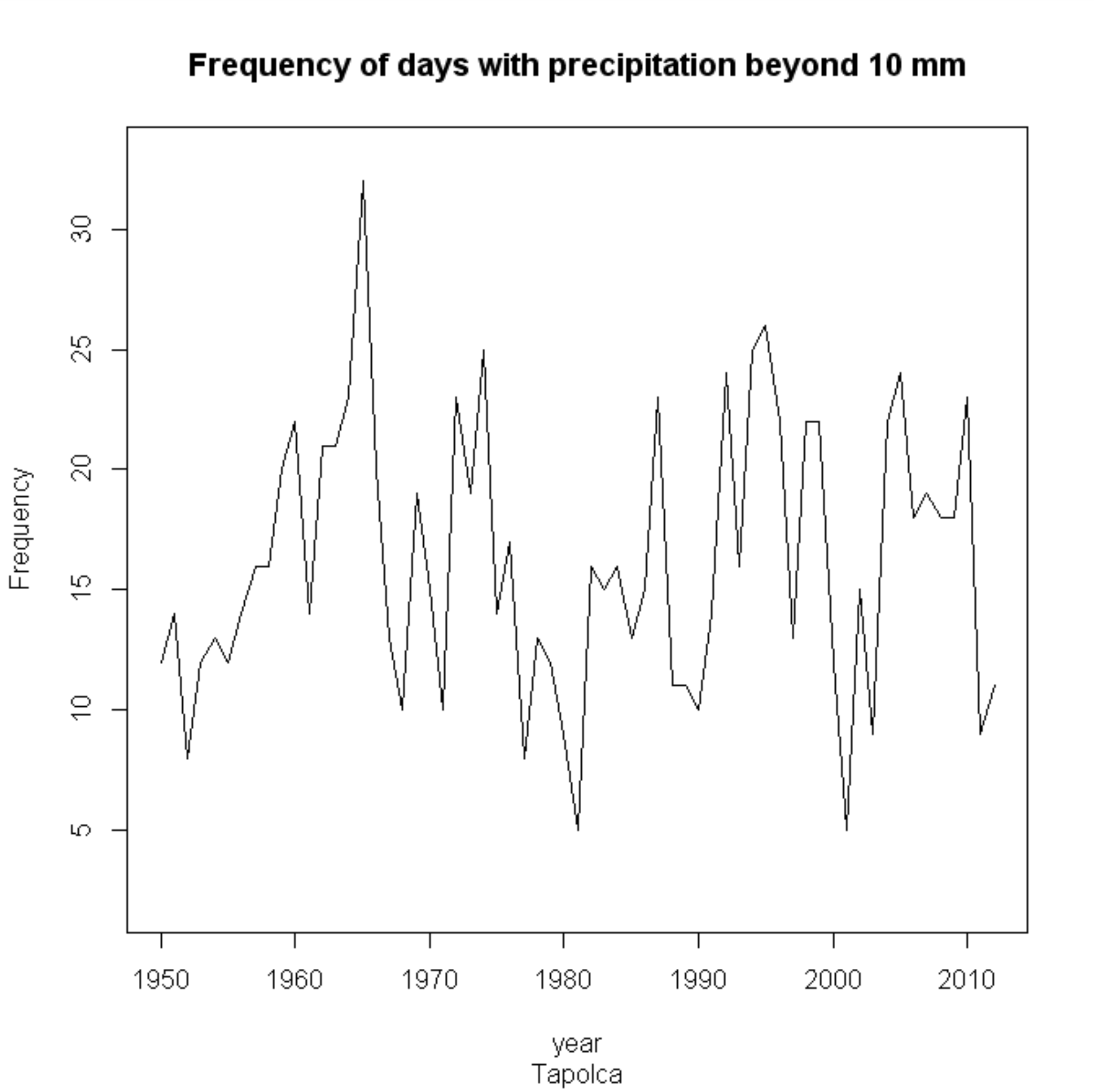}
\caption{The annual frequency of the extreme events (annual number of days with precipitation exceeding 10 mm) at one of the stations }
\label{gyak}
\end{figure}

To check if there is any tendency in the precipitation patterns, one may use two different approaches. The first deals with the frequency of the extreme events. We have carried out such an analysis, the result is shown in Figure \ref{gyak}. There is no obvious tendency beyond the random fluctuation - the other stations gave very similar results. The other approach deals with the magnitude of the extreme events. This is our choice, and it turns out to be a reasonable decision, as shown in the next subsection. 

\subsection {Univariate GPD fit}
\label{appuni}

\begin{figure}
\includegraphics[scale=0.35]{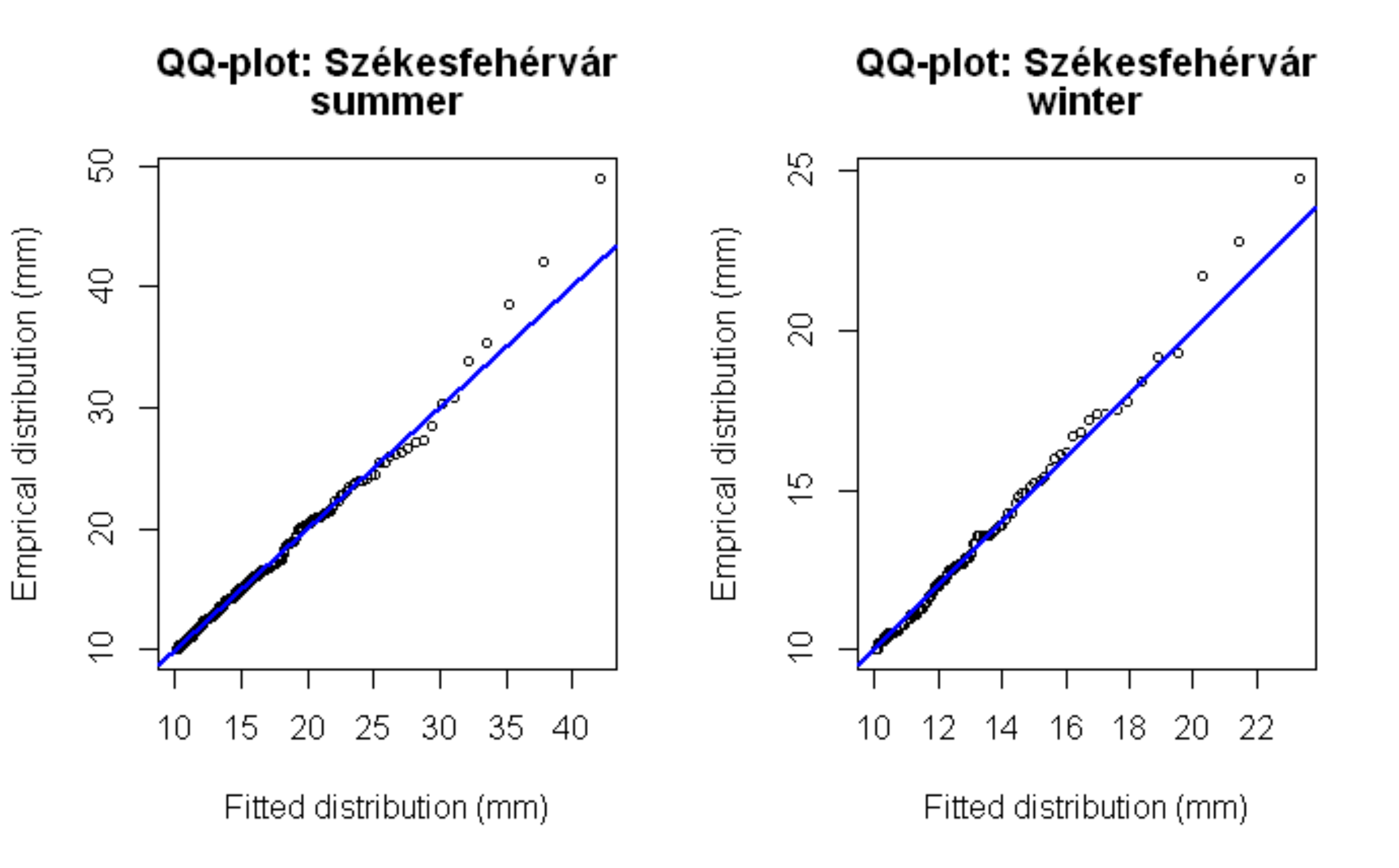}
\caption{Quantile-Quantile (QQ) plot of the GPD-fit for one of the stations}
\label{qq}
\end{figure}

We fitted the GPD separately to the summer and to the winter months. The threshold was chosen as 10 mm for simplicity in all the cases. The fit turned out to be quite good, as  Figure \ref{qq}  shows. The $y$ axis of these Quantile-Quantile plots corresponds to the observations, while the $x$ axis shows the theoretical quantiles of the fitted GPD. If the fit is good then these points should lie near to the blue line, which corresponds to the ideal case, where $y=x$.  

\begin{figure}
\includegraphics[height=80mm,angle=-90]{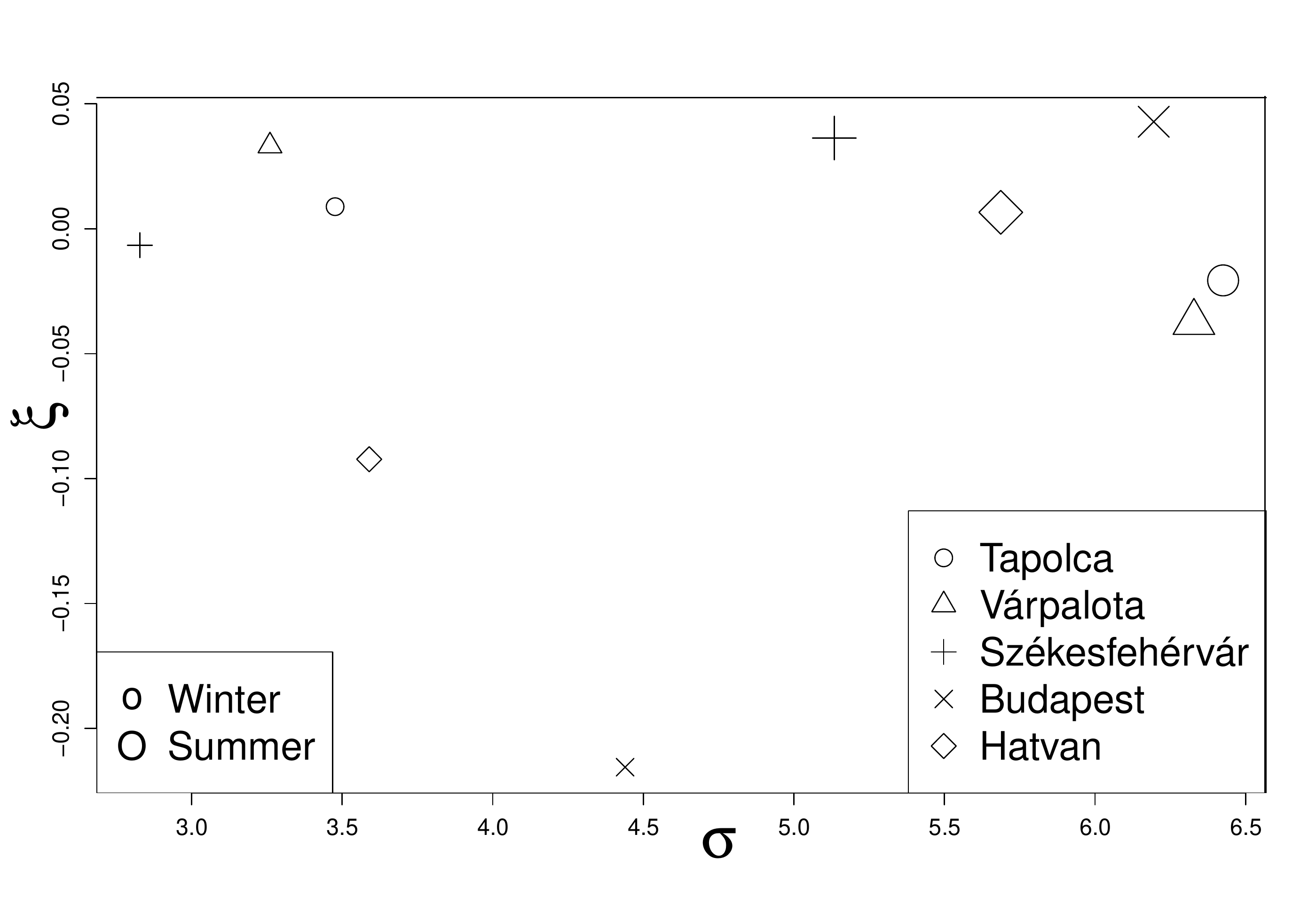}
\caption{Fitted parameters of the GPD model (see \ref{for:gpd}) for the summer and winter data, $u$=10 mm}
\label{fig:pars}
\end{figure}
In Figure \ref{fig:pars}, where the parameters of the seasons are shown,  we clearly see the differences: the scale parameter $\sigma$ is consistently larger for the summer months, and in quite a few cases also the shape parameter turns out to be higher as well.

\begin{figure*}
\includegraphics[height=88mm,angle=-90]{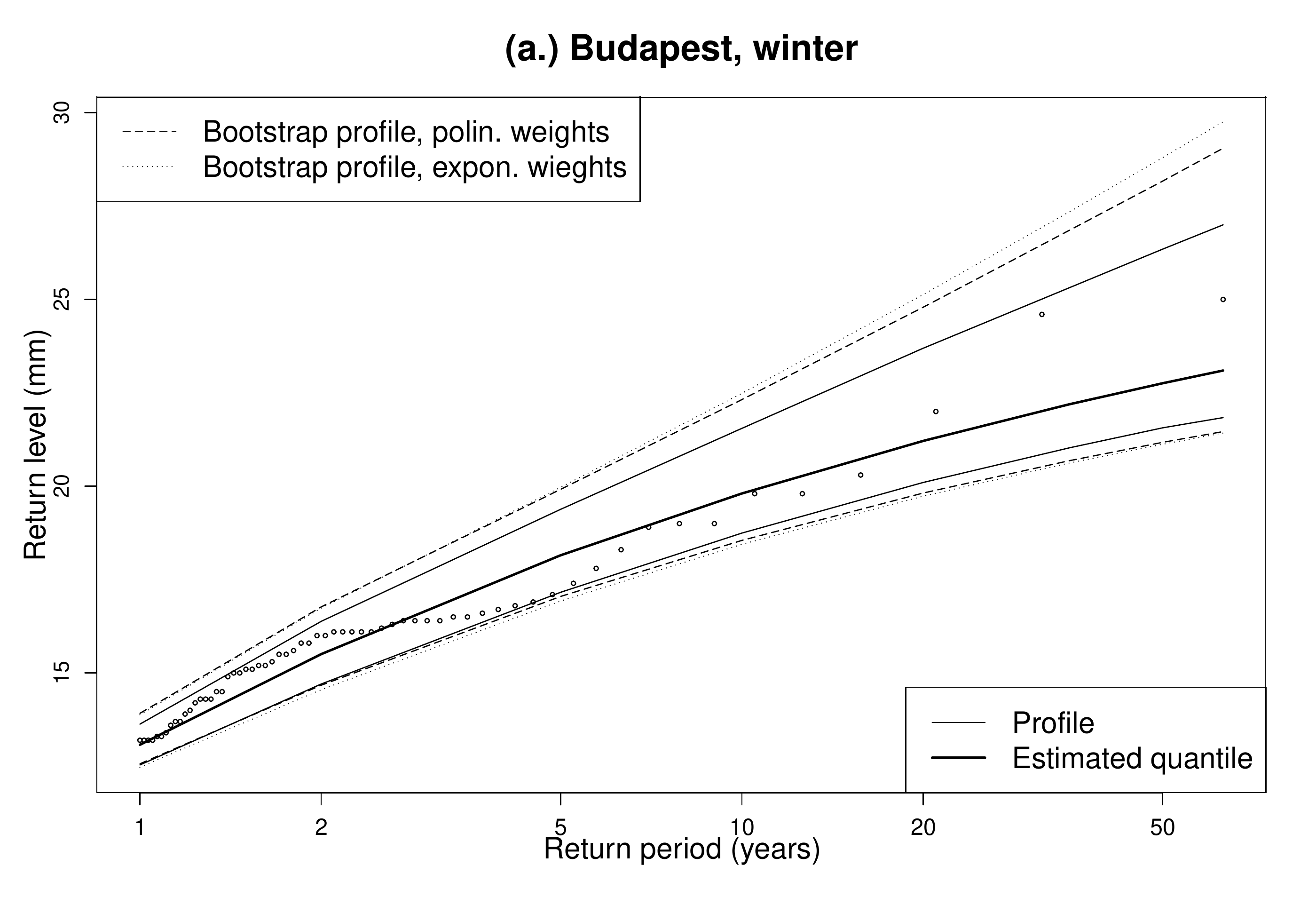}
\includegraphics[height=88mm,angle=-90]{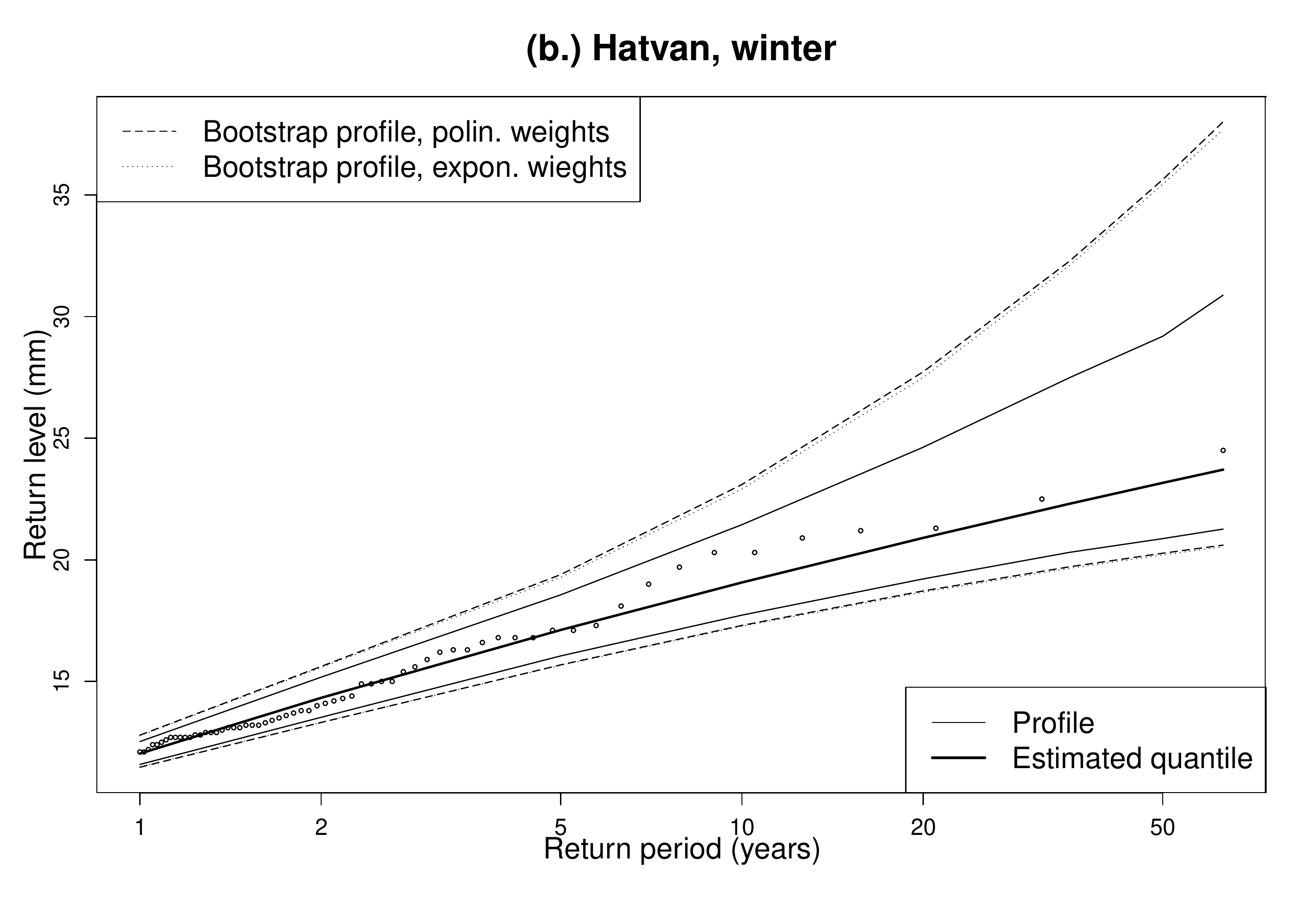} \\
\includegraphics[height=88mm,angle=-90]{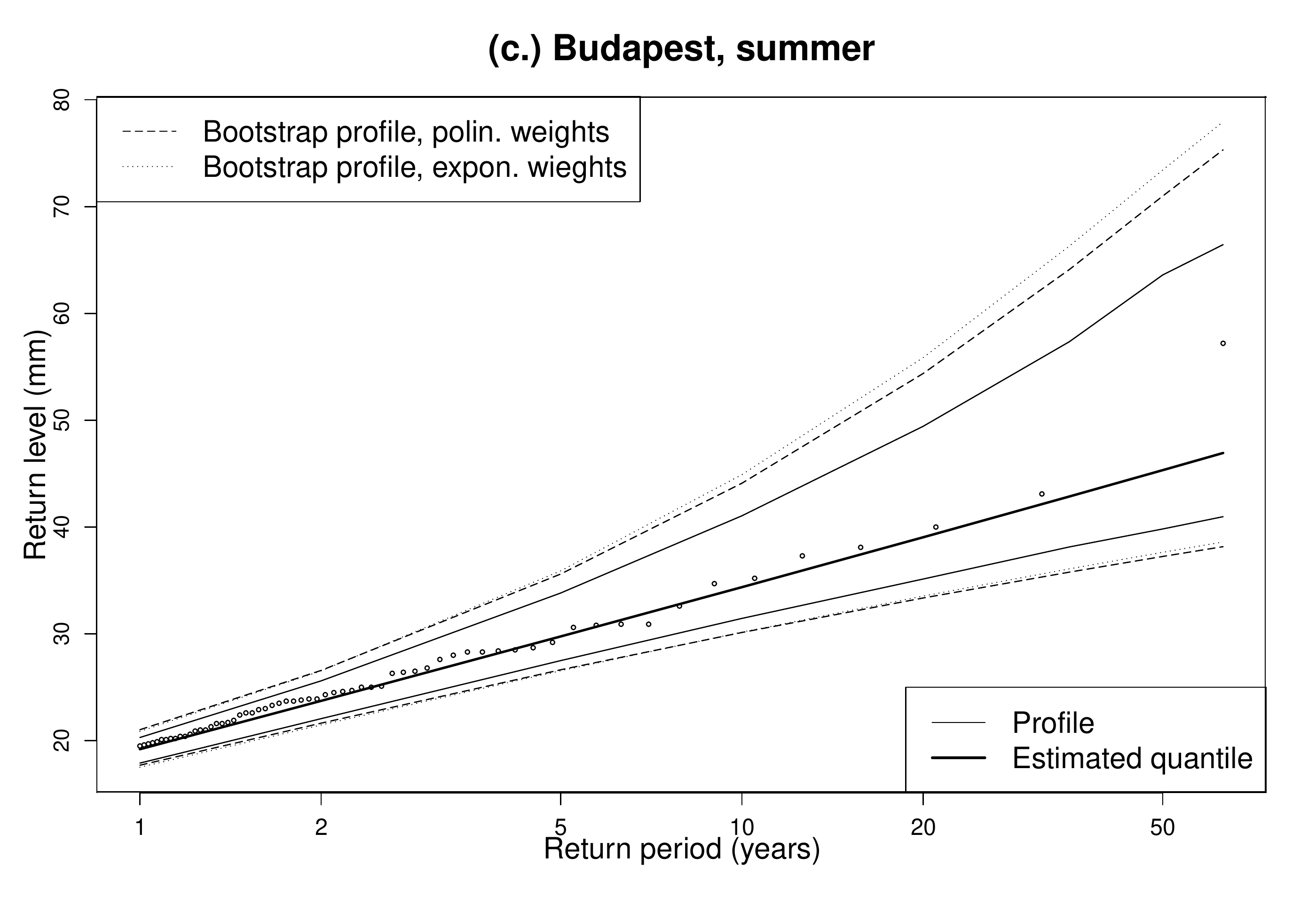}
\includegraphics[height=88mm,angle=-90]{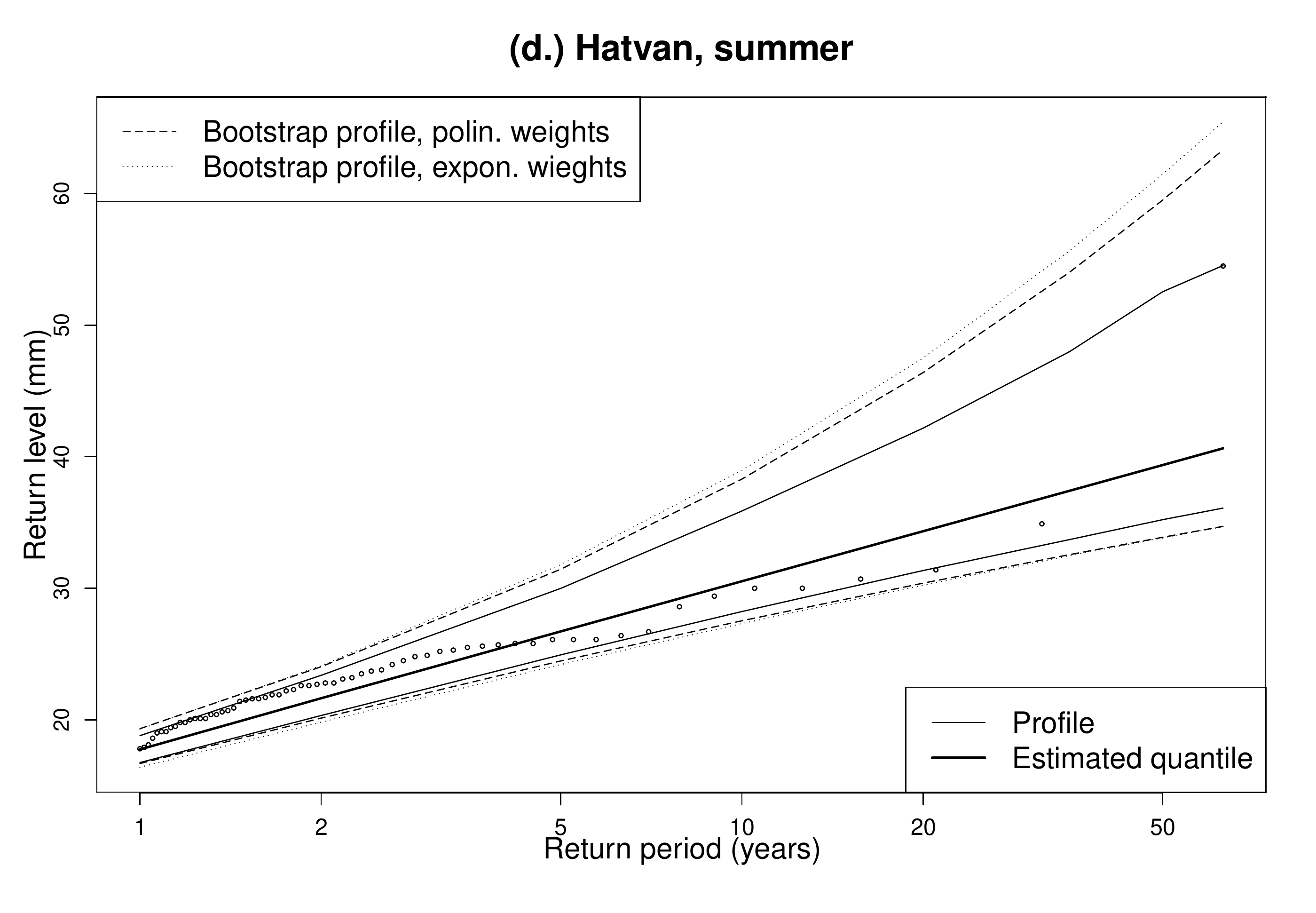} 
\vspace{-2mm}
\caption{Estimated return levels for two cities with 95\% bootstrap confidence bounds for (a.) Budapest winter data; (b.) Budapest summer data; (c.) Hatvan winter data; (d.) Hatvan summer data. The empty circles are the observations.}
\label{fig:visszat}
\end{figure*}

Next we show  results, where we present return levels together with their respective confidence bounds, based on the usual profile likelihood and our bootstrap methodology of Section \ref{boot} as well.
Figure \ref{fig:visszat} displays the return levels for the precipitation
data of Budapest and Hatvan, separately for winter (DJF) and summer (JJA). The central  line depicts the quantiles of the
estimated GPD distribution with threshold $u=10$ mm. 
The other lines show the 95\%
confidence regions by two different methods: profile likelihood and bootstrap profile likelihood.
The latter is calculated as the mean 
of the profile likelihood bounds for the weighted bootstrap samples, using exponential and multinomial
weights with the same mean and variance, respectively. We have 
used $10^4$ bootstrap replications. There does not seem to be a substantial 
difference between these choices. 

As we were interested not only in the actual results but in the properties of the methods as well, we have run a simulation study for the GPD with the actual parameter sets. 

We computed the lower and the upper 2.5\% quantiles of 100 simulated samples.

\begin{table*}[ht]
\centering
\caption{Coverage for confidence intervals, mixed GPD distributions: $w$ is the weight of the GPD with the estimated parameters for winter in Budapest, the GPD with the estimated parameters for summer has weight $1-w$}
\label{tab:sim}     
\begin{tabular}{ll|ccc}
\hline
sample size & weight ($w$) & profile lik &  weighted boot (exp)  & weighted boot (multinom) \\ \hline
50& 1& 93.3& 98.8 & 99.0 \\
100 & 1 & 94.5 & 98.7 &98.7 \\
200 & 1 &  95.2 & 98.5 & 98.6 \\
200 &     0.5&      89.8   &   97.0 &   97.4 \\ 
500 &     0.5&       88.5   &   96.7 &   97.2 \\ 
500 &     0.8   &  80.8   &  93.4  &    93.7 \\ \hline    
\end{tabular}
\end{table*}

Table \ref{tab:sim} gives the coverage for the estimated 95\% confidence intervals for return periods of 100 years. The first column shows the coverage  percentage for the classical profile likelihood, while the last two columns give the same percentiles for the weighted bootstrap methods. The results confirmed our conjecture that the often applied method of profile likelihood intervals had a low coverage if the data came from mixture distributions, like in our investigation of moving windows (subsection \ref{mw}), while our bootstrap profile methods turned out to be more conservative in these cases.

\subsection{Time dependence of the parameters}\label{mw}

\begin{figure}
\includegraphics[height=85mm,angle=-90]{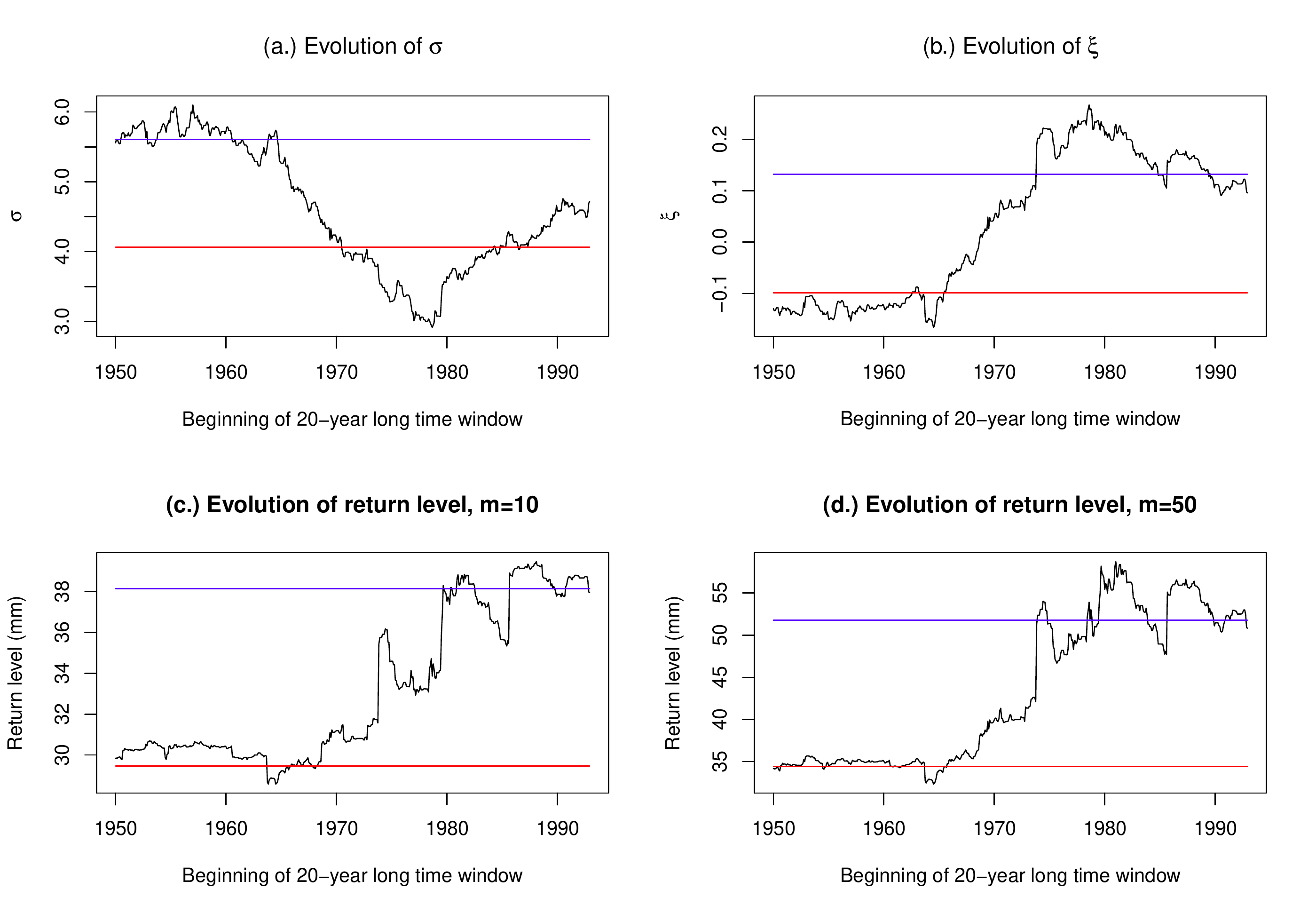} 
\caption{Results for the moving windows of length 20 years (Budapest). The blue and red lines are the upper and lower limits of the 99\% confidence intervals for the parameters; $m$ refers to the return period in years}
\label{uniwint}
\end{figure}

Our main aim is to detect if there are any significant changes in the time series. We have carried out an analysis, based on moving windows of 20 years (here we used all observations). The length of this window was chosen so that there were enough (at least 50) exceedances in every window in order to allow for reliable estimation of the parameters. We see in Figure \ref{uniwint} that there was a significant downward trend in the estimated return levels until around 1965, which was reversed afterwards, and in some cases was followed by even a significant increase. The significance level was chosen as high as 99\%, in order to reduce the chances of a type I error.  The increase in the more extreme events (50-years return levels, bottom right panel) is even more prominent than the same observation for the 10-years return periods (bottom left panel), see subsection \ref{Ch21}.

This finding seems to be new, as most of the previous authors  
\begin{itemize}
\item either investigated the annual sum of precipitation (see for example \citet{do}) and here quite often a decrease rather than an increase was observed. Our  results do not contradict to these, as we have investigated the extreme events only. The reasons behind this phenomenon together with its extent in space and time is still to be cleared; 
\item or gave forecasts with the vision of more extreme precipitation events, but it looks to us that this period with more extreme observations has already started. 
\end{itemize}


\subsection {Bivariate applications}
\label{appbi}

\begin{figure}
\includegraphics[height=82mm,angle=-90]{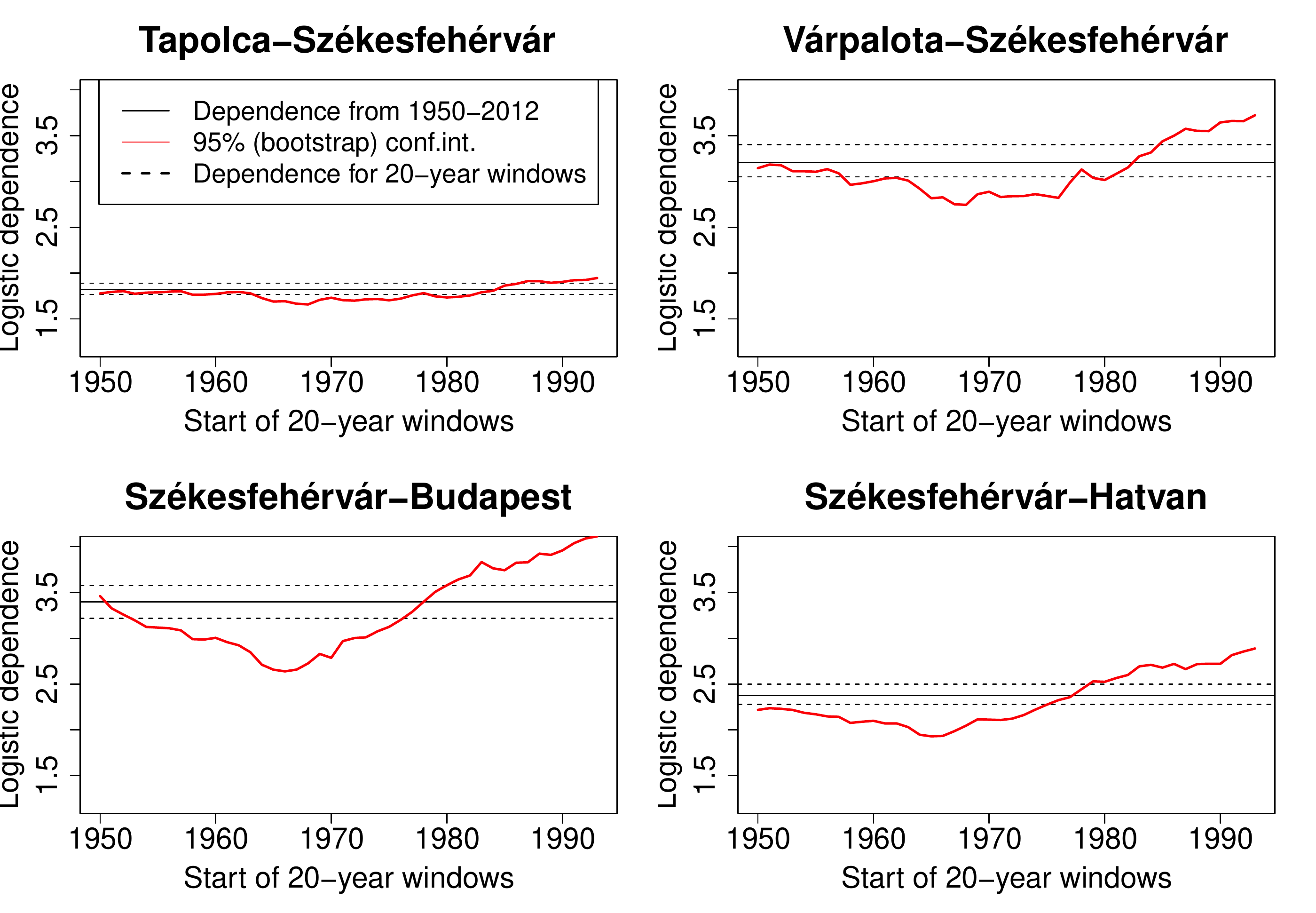}
\caption{Time dependence of the dependence parameter of the fitted logistic (BGPD II) model to four pairs}
\label{fig:mw_fvar}
\end{figure}

Seeing the results from the previous section, we now focus on the time dependence of the parameters, fitted by the R package {\tt{mgpd}} \citet{mg}.
The time dependence  was investigated by the same moving window-method as in the univariate case. Interestingly the results seen in Figure \ref{fig:mw_fvar} were also similar: in most of the cases first there was a significant decrease in the dependence, with minimal value for the period 1965-84 and afterwards an increased dependence was observed, till the most recent window of years 1993-2012. 
Table \ref{tab:mw_2d} gives the exact values of the parameters for all the pairs. We see that for all cases, except those
pairs which include Tapolca (our westernmost site) the differences are significant.

\begin{table*}[ht]
\centering
\caption{Estimated dependence parameter of the logistic BGPD II model for the entire dataset (with bootstrap confidence bounds) and for two time intervals}
\label{tab:mw_2d}     
\begin{tabular}{lccccc}
\hline
	& 	5\%	&Mean	&95\%	&1965-1984	&1993-2012\\ \hline
Tapolca-V\'arpalota	&	2.47	&2.57	&2.66	&2.43	&2.69\\
Tapolca-Sz\'ekesfeh\'erv\'ar	&1.75	&1.82	&1.86	&1.69	&1.95\\
Tapolca-Budapest	&1.43	&1.52	&1.61	&1.41	&1.59\\
Tapolca-Hatvan	&1.44	&1.54	&1.6	&1.35	&1.5\\
V\'arpalota-Sz\'ekesfeh\'erv\'ar	&3.04	&3.2	&3.34	&2.82	&3.72\\
V\'arpalota-Budapest	&1.99	&2.05	&2.11&	1.74	&2.38\\
V\'arpalota-Hatvan	&1.68	&1.73	&1.79	&1.45	&2.07\\
Sz\'ekesfeh\'erv\'ar-Budapest	&3.28	&1.82&	1.86	&2.66	&4.11\\
Sz\'ekesfeh\'erv\'ar-Hatvan	&2.29	&2.37	&2.46	&1.92	&2.89\\
Budapest-Hatvan	&3.14	&3.29	&3.48	&2.44	&3.9\\
\hline
\end{tabular}
\end{table*}

\begin{figure}
\includegraphics[height=82mm,angle=-90]{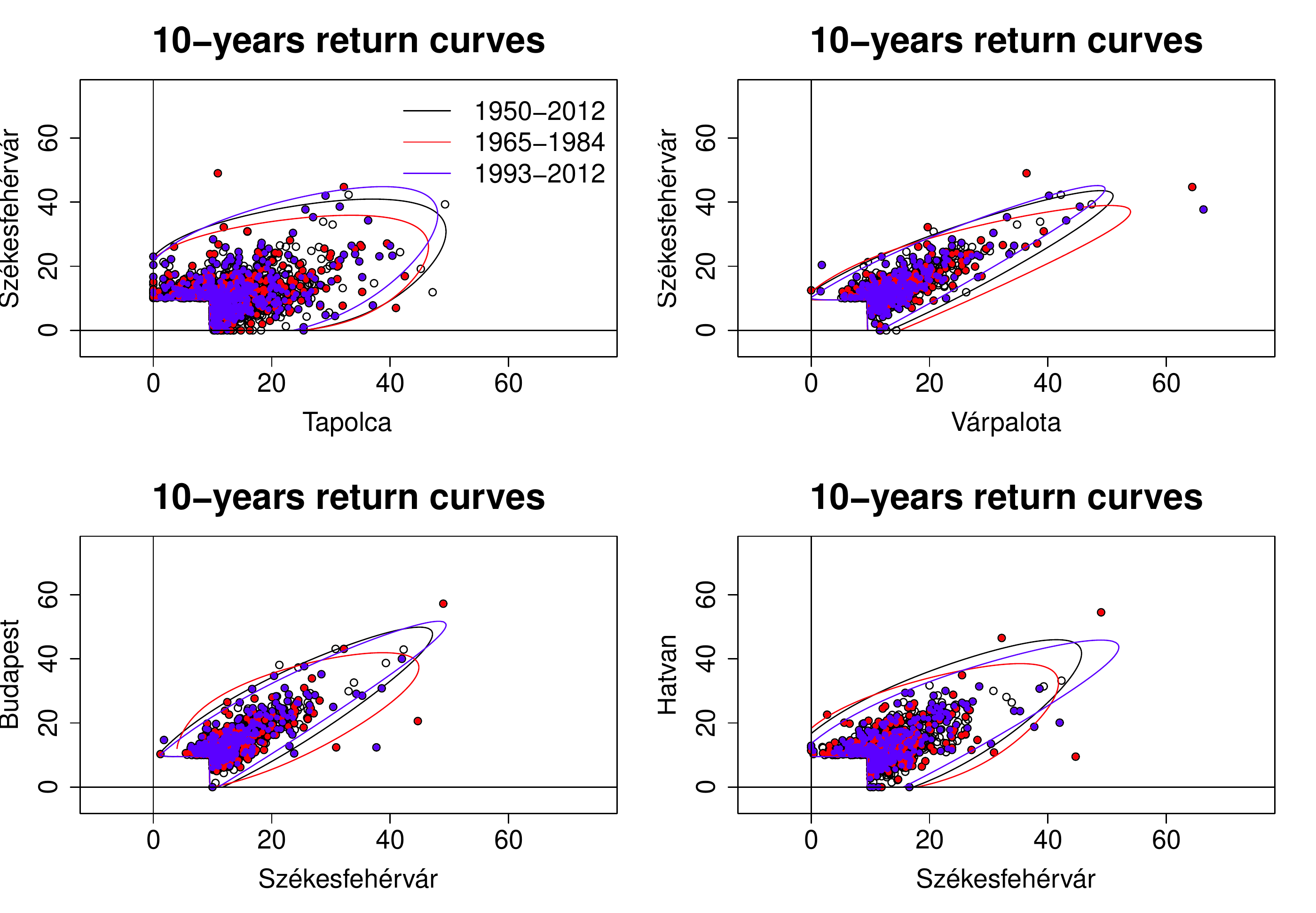}
\caption{Estimated 10 years coverage regions to four pairs; The red and blue circles correspond to the observations in the periods 1965-1984 and 1993-2012, respectively. The empty circles are the other observations.}
\label{reg}
\end{figure}

Figure \ref{reg} shows the effect of the above results to the shape of the bivariate distributions. These coverage regions were introduced in \citet{rat} as such regions which cover the given portion  of the estimated bivariate distribution and which have the smallest area among these (in our case constructed on a way that we expect an observation outside of it exactly once in 10 years). Here the changes in the estimated marginal distribution also play an important role. The differences are in some cases quite substantial. 

In Table \ref{bc}
the  bivariate quantiles for the two periods from above are compared. The values in the table give the ratio of the probabilities $p=P(X>x_q,Y>y_q)$ for the marginal $q$-quantiles for $q=0.9$. 
The increase was highly significant in quite a few cases, showing that the dependence was indeed stronger in these cases (as higher probabilities for the joint exceedances mean stronger dependence).

\begin{table*}[ht]
\centering
\caption{Estimated ratio of the probabilities of the joint exceedances of the marginal 10-year return levels in the time interval 1993-2012, in comparison to 1965-1984}
\label{bc}     
\begin{tabular}{lcccc}
\hline
	&	V\'arpalota &Sz\'ekesfeh\'erv\'ar 	&Budapest&Hatvan \\ \hline
Tapolca  & 0.631 & 0.768 & 1.020 & 0.950\\
V\'arpalota  & -- & 2.465 & 2.042 & 1.644\\
Sz\'ekesfeh\'erv\'ar     &  -- &  -- & 5.677 &  6.515\\
Budapest &   -- &   -- &   -- &11.402\\
\hline
\end{tabular}
\end{table*}
 
\section{Conclusions}
\label{conclu}

As a conclusion we can formulate that the weighted bootstrap is indeed a useful method for estimating the uncertainty even in extreme-value models.  In this case the like\-li\-hood-based approach does not cause additional computational difficulties. Another important observation is that the distributional properties of the weights beyond the second moment did not play an important role.

The BGPD II model turned out to be valuable especially in cases like the investigated one, where the number of observations is limited. The strong correlation between the mar\-gi\-nal return levels and the bivariate dependence is definitely a result which may turn out as an interesting general phenomenon.   

The findings: recently increased return levels, combined with similar observations on the increased dependence between the sites show the danger of floods -- something which has indeed been observed in summer 2013 over the Danube basin -- even if the heavy rain felt at other Central European locations this time, the tendencies might be similar, which would be worth investigating.
The advantage of the proposed methods was first the more exact quantification of the uncertainty in the confidence interval estimation via the weighted bootstrap profile likelihood and second the investigation of the extremal dependence between nearby sites. 
 
\begin{acknowledgements}
P. Rakonczai's research was supported by an {\color{white} i} OTKA mobility grant (OTKA registration number: MB08A 84576 PKR registration number: HUMAN-MB08-1-2011-0007).\\
The work of L. Varga was supported by the European Union Social Fund (Grant Agreement No. T\'AMOP 4.2.1/B-09/1/KMR-2010-0003).
\end{acknowledgements}

  
\bibliographystyle{plainnat}
\bibliography{meteorefs_1009}

\end{document}